\documentclass[10pt,twoside]{amsart}
\usepackage{amssymb}
\usepackage{float}
\usepackage{amsfonts}
\usepackage[centertags]{amsmath}
\usepackage{amsthm}
\usepackage{newlfont}
\usepackage{array}

\date{}

\begin{document}

\centerline{\Large{\bf Neutrosophic Metric Spaces}}

\centerline{}

\centerline{\bf {Murat Kiri\c{s}ci}* and \bf {Necip \c{S}im\c{s}ek}}

\centerline{}

\centerline{a. Department of Mathematical Education, Hasan Ali Y\"{u}cel Education Faculty,}

\centerline{ Istanbul University-Cerrahpa\c{s}a, Vefa, 34470, Fatih, Istanbul, Turkey}

\centerline{e-mail: mkirisci@hotmail.com}

\centerline{b. Department of mathematics, Faculty of Arts and Sciences,}

\centerline{ Istanbul Commerce University, Istanbul, Turkey}

\centerline{e-mail: necipsimsek@hotmail.com}

\centerline{}

\newtheorem{Theorem}{\quad Theorem}[section]

\newtheorem{Definition}[Theorem]{\quad Definition}

\newtheorem{Corollary}[Theorem]{\quad Corollary}

\newtheorem{Lemma}[Theorem]{\quad Lemma}

\newtheorem{Example}[Theorem]{\quad Example}

\newtheorem*{remark}{Remark}

\centerline{}
{\textbf{Abstract:}
In present paper, the definition of new metric space with neutrosophic numbers is given. Several topological and structural properties have been investigated. The analogues of Baire Category Theorem and Uniform Convergence Theorem are given for Neutrosophic metric spaces.
\\

\centerline{}

{\bf Subject Classification:} Primary 03E72; Secondary 54E35, 54A40, 46S40. \\

{\bf Keywords:}  Neutrosophic metric space, Baire Category Theorem, Uniform Convergence Theorem, nowhere dense, completeness, Hausdorffness.

\section{Introduction}

Fuzzy Sets (FSs) put forward by Zadeh \cite{Zadeh} has influenced deeply all the scientific fields since the publication of the paper.
It is seen that this concept, which is very important for real-life situations, had not enough solution to some problems in time.
New quests for such problems have been coming up. Atanassov \cite{Atan} initiated Intuitionistic fuzzy sets (IFSs) for such cases.
Neutrosophic set (NS) is a new version of the idea of the classical set which is defined by Smarandache \cite{Smar}.
Examples of other generalizations are FS \cite{Zadeh} interval-valued FS \cite{Turksen}, IFS \cite{Atan},  interval-valued IFS \cite{AtanGarg},
the sets paraconsistent, dialetheist, paradoxist, and tautological \cite{Smar0}, Pythagorean fuzzy sets \cite{Yager} .\\

Using the concepts Probabilistic metric space and fuzzy, fuzzy metric space (FMS) is introduced in \cite{KraMic}.
Kaleva and Seikkala \cite{KalSei} have defined the FMS as a distance between two points to be a non-negative fuzzy number.
In  \cite{GeoVee} some basic properties of FMS studied and the Baire Category Theorem for FMS proved.
Further, some properties such as separability, countability are given and Uniform Limit Theorem is proved in \cite{GeoVee2}.
Afterward, FMS has used in the applied sciences such as fixed point theory, image and signal processing, medical imaging, decision-making et al.
After defined of the IFS, it was used in all areas where FS theory was studied.
Park \cite{Park} defined IF metric space (IFMS), which is a generalization of FMSs.
Park used George and Veeramani's \cite{GeoVee} idea of applying t-norm and t-conorm to the FMS meanwhile defining IFMS and studying its basic features.\\

Bera and Mahapatra defined the neutrosophic soft linear spaces (NSLSs) \cite{BeraMah}. Later, neutrosophic soft normed linear spaces(NSNLS) has been defined by Bera and Mahapatra \cite{BeraMah1}. In \cite{BeraMah1}, neutrosophic norm, Cauchy sequence in NSNLS, convexity of NSNLS, metric in NSNLS were studied.\\

In present study, from the idea of neutrosophic sets, new metric space was defined which is called Neutrosophic metric Spaces (NMS). We investigate some properties of NMS such as open set, Hausdorff, neutrosophic bounded, compactness, completeness, nowhere dense. Also we give Baire Category Theorem and Uniform Convergence Theorem for NMSs.

% -----------------------------------------------------------

\section{Preliminaries}

Some definitions related to the fuzziness, intuitionistic fuzziness and neutrosophy are given as follows:\\

The fuzzy subset $F$ of $\mathbb{R}$ is said to be a fuzzy number(FN).
The FN is a mapping $F:\mathbb{R}\rightarrow [0,1]$  that corresponds to each real number $a$ to the degree of membership $F(a)$.

Let $F$ is a FN. Then, it is known that \cite{Kirisci2}

\begin{itemize}
  \item If $F(a_{0})=1$, for $a_{0} \in \mathbb{R}$, $F$ is said to be normal,
  \item If for each $\mu >0$, $F^{-1}\{[0,\tau+\mu)\}$
is open in the usual topology $\forall \tau\in [0,1)$, $F$ is said to be upper semi continuous, ,
  \item  The set $[F]^{\tau}=\{a\in\mathbb{R}:F(a)\geq \tau\}$, $\tau \in [0,1]$ is called $\tau-$cuts of $F$.
\end{itemize}

Choose non-empty set $F$. An IFS in $F$ is an object $U$ defined by
\begin{eqnarray*}
U=\{<a,G_{U}(a),Y_{U}(a)>: a\in F\}
\end{eqnarray*}
where $G_{U}(a):F\rightarrow [0,1]$ and
$Y_{U}(a):F\rightarrow [0,1]$ are functions for all $a\in F$ such that $0\leq G_{U}(a)+Y_{U}(a) \leq 1$ \cite{Atan}.
Let $U$ be an IFN. Then,

\begin{itemize}
  \item an IF subset of the $\mathbb{R}$,
  \item If $G_{U}(a_{0})=1$ and, $Y_{U}(a_{0})=0$ for $a_{0} \in \mathbb{R}$, normal,
  \item  If $G_{U}(\lambda a_{1}+(1-\lambda)a_{2})\geq \min(G_{U}(a_{1}), G_{U}(a_{2}))$, $\forall a_{1},a_{2}\in\mathbb{R}$ and $\lambda\in[0,1]$, then the membership function(MF) $G_{U}(a)$ is called convex,
  \item If $Y_{U}(\lambda a_{1}+(1-\lambda)a_{2})\geq \min(Y_{U}(a_{1}), Y_{U}(a_{2}))$, $\forall a_{1},a_{2}\in\mathbb{R}$ and $\lambda\in[0,1]$, then the nonmembership function(NMF)$Y_{U}(a)$ is concav,
  \item $G_{U}$ is upper semi continuous and $Y_{U}$ is lower semi continuous
  \item $supp U=cl(\{a\in F: Y_{U}(a)<1\})$ is bounded.
\end{itemize}

An IFS $U=\{<a,G_{U}(a),Y_{U}(a)>:a\in F\}$ such that $G_{U}(a)$ and $1-Y_{U}(a)$ are FNs, where $(1-Y_{U})(a)=1-Y_{U}(a)$, and $G_{U}(a)+Y_{U}(a)\leq 1$ is called an IFN.\\

Let's consider that $F$ is a space of points(objects). Denote the $G_{U}(a)$ is a truth-MF,
$B_{U}(a)$ is an indeterminacy-MF and $Y_{U}(a)$ is a falsity-MF, where $U$ is a set in $F$ with $a\in F$. Then, if we take $I=]0^{-},1^{+}[$
\begin{eqnarray*}
&&G_{U}(a): F\rightarrow I,\\
&&B_{U}(a): F\rightarrow I, \\
&&Y_{U}(a): F\rightarrow I,
\end{eqnarray*}

There is no restriction on the sum of $G_{U}(a)$, $B_{U}(a)$ and $Y_{U}(a)$. Therefore,
\begin{eqnarray*}
0^{-}\leq \sup G_{U}(a) + \sup B_{U}(a)+ \sup Y_{U}(a) \leq 3^{+}.
\end{eqnarray*}
The set $U$ which consist of with $G_{U}(a)$, $B_{U}(a)$ and $Y_{U}(a)$ in $F$ is called a neutrosophic sets(NS) and can be denoted by
\begin{eqnarray}\label{NS}
U=\{<a,(G_{U}(a),B_{U}(a), Y_{U}(a))>:a\in F, G_{U}(a), B_{U}(a), Y_{U}(a)\in ]0^{-},1^{+}[ \}
\end{eqnarray}

Clearly, NS is an enhancement of $[0,1]$ of IFSs.\\

An NS $U$ is included in another NS $V$, ($U\subseteq V$), if and only if,
\begin{eqnarray*}
&&\inf G_{U}(a) \leq \inf G_{V}(a), \quad  \sup G_{U}(a) \leq \sup G_{V}(a),\\
&&\inf B_{U}(a) \geq \inf B_{V}(a), \quad   \sup B_{U}(a) \geq \sup B_{V}(a),\\
&&\inf Y_{U}(a) \geq \inf Y_{V}(a),\quad   \sup Y_{U}(a) \geq \sup Y_{V}(a).
\end{eqnarray*}

for any $a\in F$. However, NSs are inconvenient to practice in real problems.
To cope with this inconvenient situation, Wang et al \cite{Wang} customized NS's definition and single-valued NSs (SVNSs) suggested.

To cope with this inconvenient situation, Wang et al \cite{Wang} customized NS's definition and single-valued NSs (SVNSs) suggested. Ye \cite{Ye}, described the notion of simplified NSs(SNSs), which may be characterized by three real numbers in the $[0,1]$. At the same time, the SNSs' operations may be impractical, in some cases \cite{Ye}.  Hence, the operations and comparison way between SNSs and the aggregation operators for SNSs are redefined in \cite{Peng}.\\

According to the Ye \cite{Ye}, a simplification of an NS $U$, in (\ref{NS}), is 
\begin{eqnarray*}
U=\left\{<a,(G_{U}(a),B_{U}(a), Y_{U}(a))>:a\in F\right\},
\end{eqnarray*}
which called an SNS. Especially, if $F$ has only one element $<G_{U}(a),B_{U}(a), Y_{U}(a)>$ is said to be an SNN. Expressly, we may see SNSs as a subclass of NSs.\\

An SNS $U$ is comprised in another SNS $V$ ($U\subseteq V$), iff $G_{U}(a)\leq G_{V}(a)$, $B_{U}(a)\geq B_{V}(a)$
and $Y_{U}(a)\geq Y_{V}(a)$ for any $a\in F$. Then, the following operations are given by Ye\cite{Ye}:

\begin{eqnarray*}
U+V&=&\langle G_{U}(a)+G_{V}(a)-G_{U}(a).G_{V}(a), B_{U}(a)+B_{V}(a)-B_{U}(a).B_{V}(a), Y_{U}(a)+Y_{V}(a)-Y_{U}(a).Y_{V}(a)\rangle,\\
U.V&=&\langle G_{U}(a).G_{V}(a), B_{U}(a).B_{V}(a), Y_{U}(a).Y_{V}(a)\rangle,\\
\alpha. U&=&\langle 1-(1-G_{U}(a))^{\alpha}, 1-(1-B_{U}(a))^{\alpha}, 1-(1-Y_{U}(a))^{\alpha}\rangle \quad \quad for \quad \alpha>0,\\
U^{\alpha}&=&\langle G_{U}^{\alpha}(a), B_{U}^{\alpha}(a), Y_{U}^{\alpha}(a)\rangle  \quad \quad for \quad \alpha>0.
\end{eqnarray*}

Triangular norms (t-norms) (TN) were initiated by Menger \cite{Menger}. In the problem of computing the distance between two elements in space, Menger offered using probability distributions instead of using numbers for distance. TNs are used to generalize with the probability distribution of triangle inequality in metric space conditions. Triangular conorms (t-conorms) (TC) know as dual operations of TNs. TNs and TCs are very significant for fuzzy operations(intersections and unions).

\begin{Definition} Give
an operation $\circ : [0,1] \times [0,1] \rightarrow [0,1]$.
If the operation $\circ$ is satisfying the following conditions, then it is called that the operation $\circ$ is \emph{continuous TN}: For $s,t,u,v\in [0,1]$, 
\begin{itemize}
\item [i.] $s \circ 1 = s$
\item [ii.] If $s\leq u$ and $t \leq v$, then $s \circ t \leq u \circ v$,
\item [iii.] $\circ$ is continuous,
\item [iv.] $\circ$ is commutative and associative.
\end{itemize}
\end{Definition}

\begin{Definition} Give
an operation $\bullet : [0,1] \times [0,1] \rightarrow [0,1]$.
If the operation $\bullet$ is satisfying the following conditions, then it is called that the operation $\bullet$ is \emph{continuous TC}:
\begin{itemize}
\item [i.] $s \bullet 0 = s$,
\item [ii.] If $s\leq u$ and $t \leq v$, then $s \bullet t \leq u \bullet v$,
\item [iii.] $\bullet $ is continuous,
\item [iv.] $\bullet$ is commutative and associative.
\end{itemize}
\end{Definition}

Form above definitions, we note that if we choose $0<\varepsilon_{1}, \varepsilon_{2}<1$ for $\varepsilon_{1}>\varepsilon_{2}$, then there exist $0<\varepsilon_{3}, \varepsilon_{4}<0,1$ such that $\varepsilon_{1}\circ \varepsilon_{3}\geq \varepsilon_{2}$, \quad $\varepsilon_{1} \geq \varepsilon_{4}\bullet \varepsilon_{2}$. Further, if we choose $\varepsilon_{5} \in (0,1)$, then there exist $\varepsilon_{6}, \varepsilon_{7} \in (0,1)$ such that $\varepsilon_{6}\circ \varepsilon_{6}\geq \varepsilon_{5}$ and $\varepsilon_{7}\bullet \varepsilon_{7}\leq \varepsilon_{5}$.

% -----------------------------------------------------------

\section{Neutrosophic Metric Spaces}

\begin{Definition}
Take $F$ be an arbitrary set, $\mathcal{N}=\{<a, G(a),B(a),Y(a)> : a\in F\}$ be a NS such that $\mathcal{N}: F \times F \times \mathbb{R}^{+} \rightarrow [0,1]$.
Let $\circ$ and $\bullet$ show the continuous TN and continuous TC, respectively.
The four-tuple $(F, \mathcal{N}, \circ, \bullet)$ is called neutrosophic metric space(NMS) when the following conditions are satisfied.
$\forall a,b,c\in F$,

\begin{itemize}
\item [i.] $0 \leq G(a,b,\lambda) \leq 1$,\quad $0 \leq B(a,b,\lambda) \leq 1$, \quad $0 \leq Y(a,b,\lambda) \leq 1$ \quad $\forall \lambda \in \mathbb{R}^{+}$,
\item [ii.] $G(a,b,\lambda)+B(a,b,\lambda)+Y(a,b,\lambda)\leq 3$,   (for $\lambda \in  \mathbb{R}^{+}$),
\item [iii.] $G(a,b,\lambda)=1$ \quad (for $\lambda >0$)  if and only if $a=b$,
\item[iv.] $G(a,b,\lambda)=G(b,a,\lambda)$ \quad (for $\lambda >0$),
\item[v.] $G(a,b,\lambda)\circ G(b,c,\mu)\leq G(a,c,\lambda+\mu)$ \quad $(\forall \lambda, \mu >0)$,
\item[vi.] $G(a,b,.): [0,\infty)\rightarrow [0,1]$ is continuous,
\item[vii.] $lim_{\lambda\rightarrow \infty}G(a,b,\lambda)=1$ \quad $(\forall\lambda >0)$,
\item [viii.] $B(a,b,\lambda)=0$ \quad (for $\lambda >0$)  if and only if $a=b$,
\item[ix.] $B(a,b,\lambda)=B(b,a,\lambda)$ \quad (for $\lambda >0$),
\item[x.] $B(a,b,\lambda)\bullet B(b,c,\mu)\geq B(a,c,\lambda+\mu)$ \quad $(\forall \lambda, \mu >0)$,
\item[xi.] $B(a,b,.): [0,\infty)\rightarrow [0,1]$ is continuous,
\item[xii.] $lim_{\lambda\rightarrow \infty}B(a,b,\lambda)=0$ \quad $(\forall\lambda >0)$,
\item [xiii.] $Y(a,b,\lambda)=0$ \quad (for $\lambda >0$)  if and only if $a=b$,
\item[xiv.] $Y(a,b,\lambda)=Y(b,a,\lambda)$ \quad $(\forall\lambda >0)$,
\item[xv.] $Y(a,b,\lambda)\bullet Y(b,c,\mu)\geq Y(a,c,\lambda+\mu)$ \quad $(\forall \lambda, \mu >0)$,
\item[xvi.] $Y(a,b,.): [0,\infty)\rightarrow [0,1]$ is continuous,
\item[xvii.] $lim_{\lambda \rightarrow \infty}Y(a,b,\lambda)=0$ \quad (for $\lambda >0$),
\item[xviii.] If $\lambda \leq 0$, then $G(a,b,\lambda)=0$, $B(a,b,\lambda)=1$ and $Y(a,b,\lambda)=1$.
\end{itemize}
Then $\mathcal{N}=(G,B,Y)$ is called Neutrosophic metric(NM) on $F$.
\end{Definition}
The functions $G(a,b,\lambda), B(a,b,\lambda), Y(a,b,\lambda)$ denote the degree of nearness, the degree of neutralness and the degree of non-nearness between $a$ and $b$ with respect to $\lambda$, respectively.

\begin{Example}\label{exp:0}
Let $(F,\textbf{d})$ be a MS. Give the operations $\circ$ and $\bullet$ as default (min) TN $a\circ b= \min \lbrace a,b\rbrace$ and default(max) TC $a \bullet b = \max \lbrace a,b\rbrace$.

\begin{eqnarray*}
G(a,b, \lambda)=\frac{\lambda}{\lambda + d(a,b)}, \quad B(a,b, \lambda)=\frac{d(a,b)}{\lambda + d(a,b)} \quad Y(a,b, \lambda)=\frac{d(a,b)}{\lambda}, 
\end{eqnarray*}
$\forall a,b\in F$ and $\lambda >0$. Then, $(F, \mathcal{N}, \circ, \bullet)$ is NMS such that  $\mathcal{N}: F\times F \times \mathbb{R}^{+} \rightarrow  [0,1]$. This NMS is expressed as produced by a metric $\textbf{d}$ the NM.

\end{Example}

\begin{Example}\label{exp:1}
Choose $F$ as natural numbers set. Give the operations $\circ$ and $\bullet$ as TN $a \circ b= \max \lbrace0, a+b-1\rbrace$ and TC $a \bullet b = a+b-ab$.
$\forall a,b\in F$, \quad $\lambda >0$
\begin{eqnarray*}
G(a,b, \lambda) = \left\{ \begin{array}{ccl}
\frac{a}{b}&, & (a\leq b),\\
\frac{b}{a}&, & (b\leq a),
\end{array} \right.
\end{eqnarray*}
\begin{eqnarray*}
B(a,b, \lambda) = \left\{ \begin{array}{ccl}
\frac{b-a}{y}&, & (ax\leq b),\\
\frac{a-b}{x}&, & (b\leq a),
\end{array} \right.
\end{eqnarray*}
\begin{eqnarray*}
Y(a,b, \lambda) = \left\{ \begin{array}{ccl}
b-a&, & (a\leq b),\\
a-b&, & (b\leq a),
\end{array} \right.
\end{eqnarray*}
Then, $(F, \mathcal{N}, \circ, \bullet)$ is NMS such that $\mathcal{N}: F\times F \times \mathbb{R}^{+} \rightarrow  [0,1]$.
\end{Example}

\begin{remark}
$\mathcal{N}=\{<a, G(a),B(a),Y(a)> : a\in F\}$ defined in Example \ref{exp:0} is not a NM with TN $a \circ b= \max \lbrace0, a+b-1\rbrace$ and TC $a \bullet b = a+b-ab$.\\

It can also be said that $\mathcal{N}=\{<a, G(a),B(a),Y(a)> : a\in F\}$ defined in Example \ref{exp:1} is not a NM with TN $a\circ b= \min \lbrace a,b\rbrace$ and TC $a \bullet b = \max \lbrace a,b\rbrace$.\\
\end{remark}

\begin{Definition}\label{ob:1}
Give $(F, \mathcal{N}, \circ, \bullet)$ be a NMS, $0<\varepsilon<1$, $\lambda >0$ and $a\in F$. The set
$O(a,\varepsilon,\lambda)=\lbrace b\in F: G(a,b,\lambda)> 1-\varepsilon, \quad  B(a,b,\lambda)< \varepsilon, \quad Y(a,b,\lambda)< \varepsilon\rbrace$ is said to be the open ball (OB) (center $a$ and radius $\varepsilon$ with respect to $\lambda$).
\end{Definition}

\begin{Theorem}\label{ob:2}
Every OB $O(a,\varepsilon,\lambda)$ is an open set (OS).
\end{Theorem}

\begin{proof}
Take $O(a,\varepsilon,\lambda)$ be an OB (center $a$, radius $\varepsilon$). Choose $b \in O(a,\varepsilon,\lambda)$.
Therefore,  $G(a,b,\lambda)> 1-\varepsilon, \quad  B(a,b,\lambda)< \varepsilon, \quad Y(a,b,\lambda)< \varepsilon$. There exists $\lambda_{0}\in (0,\lambda)$ such that  $G(a,b,\lambda_{0})> 1-\varepsilon, \quad  B(a,b,\lambda_{0})< \varepsilon, \quad Y(a,b,\lambda_{0})< \varepsilon$ because of $G(a,b,\lambda)> 1-\varepsilon$. If we take $\varepsilon_{0}=G(a,b,\lambda_{0})$, then for $\varepsilon_{0}>1-\varepsilon$,  $\zeta \in (0,1)$ will exist such that $\varepsilon_{0}>1-\zeta > 1-\varepsilon$. Give $\varepsilon_{0}$ and $\zeta$ such that $\varepsilon_{0}>1-\zeta$. Then, $\varepsilon_{1}, \varepsilon_{2}, \varepsilon_{3} \in (0,1)$ will exist
such that $\varepsilon_{0} \circ \varepsilon_{1} > 1-\zeta$, $(1-\varepsilon_{0})\bullet (1-\varepsilon_{2})\leq \zeta$ and $(1-\varepsilon_{0})\bullet (1-\varepsilon_{3})\leq \zeta$. Choose $\varepsilon_{4}=max\lbrace\varepsilon_{1}, \varepsilon_{2}, \varepsilon_{3}\rbrace$. Consider the OB $O(b,1-\varepsilon_{4},\lambda-\lambda_{0})$. We will show that $O(b,1-\varepsilon_{4},\lambda-\lambda_{0})\subset O(a,\varepsilon,\lambda)$. If we take $c \in O(b,1-\varepsilon_{4},\lambda-\lambda_{0})$, then $G(b,c,\lambda-\lambda_{0})> \varepsilon_{4}$, $B(b,c,\lambda-\lambda_{0})< \varepsilon_{4}$ and $Y(b,c,\lambda-\lambda_{0})< \varepsilon_{4}$. Then,
\begin{eqnarray*}
&& G(a,c,\lambda) \geq G(a,b,\lambda_{0}) \circ G(b,c,\lambda-\lambda_{0})\geq \varepsilon_{0} \circ \varepsilon_{4} \geq \varepsilon_{0} \circ \varepsilon_{1} \geq  1- \zeta > 1-\varepsilon,\\
&& B(a,c,\lambda) \leq B(a,b,\lambda_{0})\bullet B(b,c,\lambda-\lambda_{0}) \leq (1-\varepsilon_{0})\bullet (1-\varepsilon_{4}) \leq (1-\varepsilon_{0})\bullet (1-\varepsilon_{2}) \leq  \zeta < \varepsilon, \\
&& Y(a,c,\lambda) \leq Y(a,b,\lambda_{0})\bullet Y(b,c,\lambda-\lambda_{0}) \leq (1-\varepsilon_{0})\bullet (1-\varepsilon_{4}) \leq (1-\varepsilon_{0})\bullet (1-\varepsilon_{2}) \leq  \zeta < \varepsilon
\end{eqnarray*}

It shows that $c \in O(a,\varepsilon,\lambda)$ and $O(b,1-\varepsilon_{4},\lambda-\lambda_{0})\subset O(a,\varepsilon,\lambda)$.
\end{proof}

\begin{remark}
From the Definition \ref{ob:1} and Theorem \ref{ob:2}, we can say that $\tau_{\mathcal{N}}=\lbrace A \subset F: \textit{there exist} \quad \lambda >0 \quad \textit{and} \quad \varepsilon \in (0,1) \quad \textit{such that} \quad O(a,b,\lambda)\subset A, \quad \textit{for each} \quad a\in A \rbrace$ is a topology on $F$. In that case, every NM $\mathcal{N}$ on $F$ produces a topology $\tau_{\mathcal{N}}$ on $F$ which has a base the family of OSs of $\lbrace O(a,\varepsilon,\lambda): a\in F,  \varepsilon \in (0,1), \lambda >0 \rbrace$. This can be proved in a similar to the proof of Theorem 28 in \cite{Kirisci4}.
\end{remark}

\begin{Theorem}\label{haus:1}
Every NMS is Hausdorff.
\end{Theorem}

\begin{proof}
Let $(F, \mathcal{N}, \circ, \bullet)$ be a NMS. Choose $a$ and $b$ as two distinct points in $F$. Hence,
$0 < G(a,b,\lambda) <1$, \quad $0 < B(a,b,\lambda) < 1$, \quad $0< Y(a,b,\lambda) < 1$. Take $\varepsilon_{1}=G(a,b,\lambda)$, \quad $\varepsilon_{2}=B(a,b,\lambda)$, \quad $\varepsilon_{3}=Y(a,b,\lambda)$ and $\varepsilon=max\lbrace \varepsilon_{1}, 1-\varepsilon_{2}, 1-\varepsilon_{3} \rbrace$. If we take $\varepsilon_{0}\in (\varepsilon, 1)$, then there exist $\varepsilon_{4}, \varepsilon_{5}, \varepsilon_{6}$ such that $\varepsilon_{4}\circ \varepsilon_{4} \geq \varepsilon_{0}$, $(1-\varepsilon_{5})\bullet (1-\varepsilon_{5})\leq 1-\varepsilon_{0}$ and $(1-\varepsilon_{6})\bullet (1-\varepsilon_{6})\leq 1-\varepsilon_{0}$. Let $\varepsilon_{7}=max\lbrace \varepsilon_{4}, \varepsilon_{5}, \varepsilon_{6} \rbrace$. If we consider the OBs $O(a,1-\varepsilon_{7},\frac{\lambda}{2})$ and $O(b,1-\varepsilon_{7},\frac{\lambda}{2})$, then clearly $O(a,1-\varepsilon_{7},\frac{\lambda}{2}) \bigcap O(b,1-\varepsilon_{7},\frac{\lambda}{2})= \emptyset$. From here, if we choose $c\in O(a,1-\varepsilon_{7},\frac{\lambda}{2}) \bigcap O(b,1-\varepsilon_{7},\frac{\lambda}{2})$, then

\begin{eqnarray*}
&&\varepsilon_{1}=G(a,b,\lambda)\geq G(a,c,\frac{\lambda}{2})\circ G(c,b,\frac{\lambda}{2})\geq \varepsilon_{7} \circ \varepsilon_{7} \geq \varepsilon_{4} \circ \varepsilon_{4} \geq \varepsilon_{0} > \varepsilon_{1},\\
&&\varepsilon_{2}=B(a,b,\lambda)\leq B(a,c,\frac{\lambda}{2})\bullet B(c,b,\frac{\lambda}{2})\leq (1-\varepsilon_{7})\bullet (1-\varepsilon_{7})\leq (1-\varepsilon_{5})\bullet (1-\varepsilon_{5})\leq 1-\varepsilon_{0}<\varepsilon_{2},
\end{eqnarray*}
and
\begin{eqnarray*}
&&\varepsilon_{3}=Y(a,b,\lambda)\leq Y(a,c,\frac{\lambda}{2})\bullet Y(c,b,\frac{\lambda}{2})\leq (1-\varepsilon_{7})\bullet (1-\varepsilon_{7})\leq (1-\varepsilon_{6})\bullet (1-\varepsilon_{6})\leq 1-\varepsilon_{0}<\varepsilon_{3},
\end{eqnarray*}
which is a contradiction. Therefore, we say that NMS is Hausdorff.
\end{proof}

\begin{Definition}
Let $(F, \mathcal{N}, \circ, \bullet)$ be a NMS. A subset $A$ of $F$ is called Neutrosophic-bounded (NB),
 if there exist $\lambda >0$ and $\varepsilon \in (0,1)$ such that $G(a,b,\lambda)> 1-\varepsilon$,
$B(a,b,\lambda)< \varepsilon$ and $Y(a,b,\lambda)< \varepsilon$  \quad  ($\forall a,b\in A$).
\end{Definition}

\begin{Definition}
If $A \subseteq \cup_{U\in \mathcal{C}_{\mathcal{N}}}U$, a collection $\mathcal{C}_{\mathcal{N}}$ of OSs is said to be an open cover(OC) of $A$. A subspace $A$ of a NMS is compact, if every OC of $A$ has a finite subcover.\\

If every sequence in $A$ has a convergent subsequence to a point in $A$, then it is called sequential compact.
\end{Definition}

\begin{Theorem}\label{comp:1}
Every compact subset $A$ of a NMS is NB.
\end{Theorem}

\begin{proof}
Firstly, choose a compact subset $A$ of NMS $F$. Consider the OC $\lbrace O(a, \varepsilon, \lambda): a\in A\rbrace$ for $\lambda >0$, \quad $\varepsilon \in (0,1)$. Since $A$ is compact, then there exist $a_{1}, a_{2}, \ldots, a_{n} \in A$ such that $A \subseteq \cup_{k=1}^{n}O(a_{k}, \varepsilon, \lambda)$.
For some $k,m$ and $a,b\in A$, $a\in O(a_{k}, \varepsilon, \lambda)$ and $b\in O(a_{m}, \varepsilon, \lambda)$.
Then we can write, $G(a, a_{k}, \lambda)> 1-\varepsilon$, \quad $B(a, a_{k}, \lambda)< \varepsilon$, \quad $Y(a, a_{k}, \lambda)< \varepsilon$ and $G(b, a_{m}, \lambda)> 1-\varepsilon$, \quad $B(b, a_{m}, \lambda)< \varepsilon$, \quad $Y(b, a_{m}, \lambda)< \varepsilon$. Let $\rho=min\lbrace G(a_{k},a_{m}, \lambda): 1\leq k, m\leq n \rbrace$, \quad $\sigma =max\lbrace B(a_{k},a_{m}, \lambda): 1\leq k, m\leq n \rbrace$ and $\varphi=max\lbrace Y(a_{k},a_{m}, \lambda): 1\leq k, m\leq n \rbrace$. Then, $\rho, \sigma, \varphi>0$. From here, for $0< \zeta_{1}, \zeta_{2},\zeta_{3} <1$,
\begin{eqnarray*}
&&G(a,b, 3\lambda)\geq G(a,a_{k},\lambda)\circ G(a_{k},a_{m},\lambda) \circ G(a_{m},b,\lambda)\geq (1-\varepsilon) \circ (1-\varepsilon) \circ \rho > 1- \zeta_{1},\\
&&B(a,b, 3\lambda)\leq B(a,a_{k},\lambda)\bullet B(a_{k},a_{m},\lambda) \bullet B(a_{m},b,\lambda) \leq \varepsilon \bullet \varepsilon \bullet \sigma < \zeta_{2},\\
&&Y(a,b, 3\lambda)\leq Y(a,a_{k},\lambda)\bullet Y(a_{k},a_{m},\lambda) \bullet Y(a_{m},b,\lambda) \leq \varepsilon \bullet \varepsilon \bullet \varphi < \zeta_{3}.
\end{eqnarray*}

If we take $\zeta=max\lbrace \zeta_{1}, \zeta_{2}, \zeta_{3} \rbrace$ and $\lambda_{0}=3\lambda$, we have $G(a,b, \lambda_{0})>1-\zeta$, \quad $B(a,b, \lambda_{0})<\zeta$ and $Y(a,b, \lambda_{0})<\zeta$ for all $a,b\in A$. This result leads us to the conclusion that the set $A$ is NB.
\end{proof}

If $(FX, \mathcal{N}, \circ, \bullet)$ is NMS produces by a metric $\textbf{d}$ on $X$ and $A\subset F$, then
$A$ is NB if and only if it is bounded. Consequently, with Theorems \ref{haus:1} and \ref{comp:1},
we can write:\\

\begin{Corollary}
In a NMS, every compact set is closed and bounded.
\end{Corollary}

\begin{Theorem}
Take $(F, \mathcal{N}, \circ, \bullet)$ is FMS and $\tau_{\mathcal{N}}$ be the topology on $F$ produced by the FM. Then for a sequence $(a_{n})$ in $F$, the sequence $a_{n}$ is convergent to $a$ if and only if $G(a_{n},a,\lambda) \rightarrow 1$, $B(a_{n},a,\lambda)\rightarrow 0$ and $Y(a_{n},a,\lambda)\rightarrow 0$ as $n \rightarrow \infty$.
\end{Theorem}

\begin{proof}
Take $\lambda >0$. Assume that $a_{n}\rightarrow a$. If $0<\varepsilon<1$, then there exist $N\in \mathbb{N}$ such that $a_{n} \in O(a, \varepsilon, \lambda)$, \quad $ (\forall n \geq N $). Therefore,
$ 1- G(a_{n},a,\lambda)<\varepsilon $, \quad $B(a_{n},a,\lambda)< \varepsilon$ and $Y(a_{n},a,\lambda)< \varepsilon$. In that case, we can write $G(a_{n},a,\lambda)\rightarrow 1$, \quad $B(a_{n},a,\lambda)\rightarrow 0$ and $Y(a_{n},a,\lambda)\rightarrow 0$ as $n\rightarrow \infty$.\\

Conversely, $G(a_{n},a,\lambda) \rightarrow 1$, $B(a_{n},a,\lambda)\rightarrow 0$ and $Y(a_{n},a,\lambda)\rightarrow 0$ as $n \rightarrow \infty$, for each $\lambda >0$. Then, for $0<\varepsilon<1$, there exist $N \in \mathbb{N}$ such that $ 1- G(a_{n},a,\lambda)<\varepsilon $, \quad $B(a_{n},a,\lambda)< \varepsilon$ and $Y(a_{n},a,\lambda)< \varepsilon$  $\forall N \in \textbf{N}$. Then, $ G(a_{n},a,\lambda)> 1- \varepsilon $, \quad $B(a_{n},a,\lambda)< \varepsilon$ and $Y(a_{n},a,\lambda)< \varepsilon$, \quad $\forall N \in \mathbb{N}$. Then, $a_{n} \in O(a,\varepsilon, \lambda)$ $\forall n \geq N$. This is the desired result.
\end{proof}

\begin{Definition}
Take $(F, \mathcal{N}, \circ, \bullet)$ to be a NMS. A sequence $(a_{n})$ in $F$ is called \textbf{Cauchy} if for each $\varepsilon >0$ and each $\lambda >0$, there exist $N \in \mathbb{N}$ such that $G(a_{n},a_{m},\lambda)> 1-\varepsilon$, $B(a_{n},a_{m},\lambda)<\varepsilon$, $Y(a_{n},a_{m},\lambda)<\varepsilon$ for all $n,m\geq N$. $(F, \mathcal{N}, \circ, \bullet)$ is called \textbf{complete} if every Cauchy sequence is convergent with respect to $\tau_{\mathcal{N}}$.
\end{Definition}

\begin{Theorem}
Take $(F, \mathcal{N}, \circ, \bullet)$ to be a NMS. Let's every Cauchy sequence in $F$ has a convergent subsequences.
Then the NMS $(F, \mathcal{N}, \circ, \bullet)$ is complete.
\end{Theorem}

\begin{proof}
Let the sequence $(a_{n})$ be a Cauchy and let $(a_{i_{n}})$ be a subsequence of $(a_{n})$ and $a_{n} \rightarrow a$. Let $\lambda >0$ and $\mu \in (0,1)$. Take $0<\varepsilon <1$ such that $(1-\varepsilon) \circ (1-\varepsilon)\geq 1-\mu$, \quad $\varepsilon \bullet \varepsilon \leq \mu$. It is known that the sequence $(a_{n})$ is Cauchy. Then there is $N \in \mathbb{N}$ such that $G(a_{m}, a_{n}, \frac{\lambda}{2})>1-\varepsilon$, $B(a_{m}, a_{n}, \frac{\lambda}{2})<\varepsilon$ and $Y(a_{m}, a_{n}, \frac{\lambda}{2})<\varepsilon$ \quad $\forall m,n\in N$. Since $a_{n_{i}}\rightarrow a$, there is positive integer $i_{p}$ such that $i_{p}>N$, $G(a_{i_{p}}, a, \frac{\lambda}{2})>1-\varepsilon$, $B(a_{i_{p}}, a, \frac{\lambda}{2})<\varepsilon$ and $Y(a_{i_{p}}, a, \frac{\lambda}{2})<\varepsilon$. Therefore, if $n\geq N$,
\begin{eqnarray*}
&&G(a_{n}, a, \lambda) \geq G(a_{n}, a_{i_{p}}, \frac{\lambda}{2})\circ G(a_{i_{p}}, a, \frac{\lambda}{3})> (1-\varepsilon) \circ (1-\varepsilon)\geq 1-\mu,\\
&&B(a_{n}, a, \lambda)\leq B(a_{n}, a_{i_{p}}, \frac{\lambda}{2})\bullet B(a_{i_{p}, a}, \frac{\lambda}{3})<\varepsilon \bullet \varepsilon \leq \mu,\\
&&Y(a_{n}, a, \lambda)\leq Y(a_{n}, a_{i_{p}}, \frac{\lambda}{2})\bullet Y(a_{i_{p}, a}, \frac{\lambda}{3})<\varepsilon \bullet \varepsilon \leq \mu.
\end{eqnarray*}
Thus, we have $a_{n} \rightarrow a$. This is the desired result.
\end{proof}

\begin{Theorem}
Let $(F, \mathcal{N}, \circ, \bullet)$ is NMS and let $A$ be a subset of $F$ with the subspace NM \quad $(G_{A}, B_{A}, Y_{A})=(G|_{A^{2}\times (0,\infty)}, B|_{A^{2}\times (0,\infty)}, Y|_{A^{2}\times (0,\infty)})$. Then
$(A, \mathcal{N}_{A}, \circ, \bullet)$ is complete if and only if $A$ is closed subset of $F$.
\end{Theorem}

\begin{proof}
Assume that $A$ is a closed subset of $F$. Choose the sequence $(a_{n})$ be a Cauchy in $(A, \mathcal{N}_{A}, \circ, \bullet)$. Since $(a_{n})$ is a Cauchy in $F$, then there is a point $a$  in $F$ such that $a_{n}\rightarrow a$. From here, $a \in \overline{A}=A$ and so $(a_{n})$ converges to $A$.\\

Contrarily, consider the $(A, \mathcal{N}_{A}, \circ, \bullet)$ is complete. Further, assume that $A$ is not closed. Choose $a\in \overline{A} / A$. Therefore, there exist a sequence $(a_{n})$ of points in $A$ that converges to $a$ and so $(a_{n})$ is a Cauchy. Hence, for $n,m \geq N$, each $0<\mu <1$, each $\lambda>0$, there is $N\in \textbf{N}$ such that $G(a_{n}, a_{m}, \lambda)>1-\mu$, $B(a_{n}, a_{m}, \lambda)<\mu$ and $Y(a_{n}, a_{m}, \lambda)<\mu$. Now, we can write $G(a_{n}, a_{m}, \lambda)=G_{A}(a_{n}, a_{m}, \lambda)$, $B(a_{n}, a_{m}, \lambda)=B_{A}(a_{n}, a_{m}, \lambda)$ and $Y(a_{n}, a_{m}, \lambda)=Y_{A}(a_{n}, a_{m}, \lambda)$ because of the sequence $(a_{n})$ is in $A$. Therefore $(a_{n})$ is a Cauchy in $A$. Since we know that  $(F, \mathcal{N}, \circ, \bullet)$ is complete, then there is a $b\in A$ such that $a_{n}\rightarrow b$. Hence, there is $N\in \textbf{N}$ such that $G_{A}(b, a_{n},\lambda)>1-\mu$, $B_{A}(b, a_{n},\lambda)<\mu$ and $Y_{A}(b, a_{n},\lambda)<\mu$ for $n \geq N$, each $0<\mu <1$ and each $\lambda>0$. Since the sequence $(a_{n})$ is in $A$ and $b\in A$, we can write $G(b, a_{n},\lambda)=G_{A}(b, a_{n},\lambda)$, $B(b, a_{n},\lambda)=B_{A}(b, a_{n},\lambda)$ and $Y(b, a_{n},\lambda)=Y_{A}(b, a_{n},\lambda)$. This gives us the conclusion that the sequence $(a_{n})$ converges to both $a$ and $b$ in $(F, \mathcal{N}, \circ, \bullet)$. Since $a \not \in A$ and $b\in A$, we have $a\neq b$. This is a contradiction and thus the desired result is achieved.
\end{proof}

In proof of Lemma \ref{lem:1} and Theorem \ref{teo:1}, used similar proof techniques of Propositions 4.3 and 4.4 in \cite{Kirisci3}.

\begin{Lemma}\label{lem:1}
Let $(F, \mathcal{N}, \circ, \bullet)$ is NMS. If $\lambda >0$ and $\varepsilon_{1}, \varepsilon_{2} \in (0,1)$ such that $(1-\varepsilon_{2})\circ (1-\varepsilon_{2})\geq (1-\varepsilon_{1})$ and $\varepsilon_{2} \bullet \varepsilon_{2}\leq \varepsilon_{1}$, then $\overline{O(a,\varepsilon_{2},\frac{\lambda}{2})}\subset O(a,\varepsilon_{1},\lambda)$.
\end{Lemma}

\begin{proof}
Let $b\in \overline{O(a,\varepsilon_{2},\frac{\lambda}{2})}$ and let $O(b,\varepsilon_{2},\frac{\lambda}{2})$ be an OB with center $a$ and radius $\varepsilon_{2}$.
Since $O(b,\varepsilon_{2},\frac{\lambda}{2})\cap O(a,\varepsilon_{2},\frac{\lambda}{2})\neq \emptyset$, there is a $c\in O(b,\varepsilon_{2},\frac{\lambda}{2})\cap O(a,\varepsilon_{2},\frac{\lambda}{2})$.
Then, we obtain
\begin{eqnarray*}
&&G(a,b,\lambda) \geq G(a,c, \frac{\lambda}{2})\circ G(b,c, \frac{\lambda}{2}) > (1-\varepsilon_{2})\circ (1-\varepsilon_{2})\geq 1-\varepsilon_{1},\\
&&B(a,b,\lambda) \leq B(a,c, \frac{\lambda}{2})\bullet B(b,c, \frac{\lambda}{2}) < \varepsilon_{2} \bullet \varepsilon_{2} \leq \varepsilon_{1},\\
&&Y(a,b,\lambda) \leq Y(a,c, \frac{\lambda}{2})\bullet Y(b,c, \frac{\lambda}{2}) < \varepsilon_{2} \bullet \varepsilon_{2} \leq \varepsilon_{1}.
\end{eqnarray*}
Hence, $c\in O(a,\varepsilon_{1}, \lambda)$ and thus $\overline{O(a,\varepsilon_{2},\frac{\lambda}{2})}\subset O(a,\varepsilon_{1},\lambda)$.
\end{proof}

\begin{Theorem}\label{teo:1}
A subset $A$ of a NMS $(F, \mathcal{N}, \circ, \bullet)$ is nowhere dense if and only if every nonempty OS in $F$ includes an OB whose closure is disjoint from $A$.
\end{Theorem}

\begin{proof}
Let $\gamma$ be a nonempty open subset of $F$. Then there exist a nonempty OS $\delta$ such that $\delta \subset \gamma$, \quad $\delta \cap \overline{A} \neq \emptyset$. If we take $a \in \delta$, then there exist $\varepsilon_{1} \in (0,1)$, \quad $\lambda >0$ such that $O(a, \varepsilon_{1}, \lambda) \subset \delta$. Now we take $\varepsilon_{2} \in (0,1)$ such that $(1-\varepsilon_{2})\circ (1-\varepsilon_{2})\geq 1-\varepsilon_{1}$ and $\varepsilon_{2} \bullet \varepsilon_{2} \leq \varepsilon_{1}$. Using the Lemma \ref{lem:1}, we have $\overline{O(a,\varepsilon_{2},\frac{\lambda}{2})}\subset O(a,\varepsilon_{1},\lambda)$. In that case, we can write $O(a,\varepsilon_{2},\frac{\lambda}{2}) \subset \gamma$ and $\overline{O(a,\varepsilon_{2},\frac{\lambda}{2})} \cap A = \emptyset$.\\

Conversely, assume that $A$ is not nowhere dense. Therefore, $int(\overline{A})\neq \emptyset$, so there exists a nonempty OS $\gamma$ such that $\gamma \subset \overline{A}$. Take $O(a, \varepsilon_{1}, \lambda)$ be an OB such that $O(a, \varepsilon_{1}, \lambda) \subset \gamma$. Then, $\overline{O(a,\varepsilon_{2},\lambda)} \cap A \neq \emptyset$.  This result indicates that there is a contradiction.
\end{proof}

Now, we will prove Baire Category Theorem for NMS:

\begin{Theorem}
Let $\{\gamma_{n}:n\in \mathbb{N}\}$ be a sequence of dense open subsets of a complete NMS $(F, \mathcal{N}, \circ, \bullet)$. Then $\cap_{n\in \mathbb{N}}\gamma_{n}$ is also dense in $F$.
\end{Theorem}

\begin{proof}
Choose $\delta$ be nonempty OS of $F$. Since $\gamma_{1}$ is dense in $F$, $\delta\cap \gamma_{1}\neq \emptyset$. Let $a_{1}\in \delta\cap \gamma_{1}$. Since $\delta\cap \gamma_{1}$ is open, then there exist $\varepsilon_{1}\in (0,1)$,\quad $\lambda_{1}>0$ such that $O(a_{1}, \varepsilon_{1}, \lambda_{1})\subset \delta \cap \gamma_{1}$. Take $\varepsilon_{1}^{*}<\varepsilon_{1}$ and $\lambda_{1}^{*}=\min\{\lambda_{1},1\}$ such that $\overline{O(a_{1}, \varepsilon_{1}^{*}, \lambda_{1}^{*})}\subset \delta \cap \gamma_{1}$. Since $\gamma_{2}$ is dense in $F$, $O(a_{1}, \varepsilon_{1}^{*}, \lambda_{1}^{*})\cap \gamma_{2}\neq \emptyset$. Let $a_{2}\in O(a_{1}, \varepsilon_{1}^{*}, \lambda_{1}^{*}) \cap \gamma_{2}$. Since $O(a_{1}, \varepsilon_{1}^{*}, \lambda_{1}^{*}) \cap \gamma_{2}$ is open, then there exist $\varepsilon_{2}\in (0, 1/2)$ and $\lambda_{2}>0$ such that $O(a_{2}, \varepsilon_{2}, \lambda_{2})\subset O(a_{1}, \varepsilon_{1}^{*}, \lambda_{1}^{*})\cap \gamma_{2}$. Take $\varepsilon_{2}^{*}<\varepsilon_{2}$ and $\lambda_{2}^{*}=\min\{\lambda_{2}, 1/2\}$ such that $\overline{O(a_{2}, \varepsilon_{2}^{*}, \lambda_{2}^{*})}\subset O(a_{1}, \varepsilon_{1}^{*}, \lambda_{2}^{*}) \cap \gamma_{2}$. If we continue this way, we have a sequence $(a_{n})$ in $F$ and a sequence $(\lambda_{n}^{*})$ such that $0<\lambda_{n}^{*}<1/n$ and
\begin{eqnarray*}
\overline{O(a_{n}, \varepsilon_{n}^{*}, \lambda_{n}^{*})}\subset O(a_{n-1}, \varepsilon_{n-1}^{*}, \lambda_{n-1}^{*})\cap \gamma_{n}
\end{eqnarray*}
Now, we show that the sequence $(a_{n})$ is a Cauchy sequence. For $\lambda >0$ and $\mu>0$, take $N\in \mathbb{N}$ such that $1/N < \lambda $ and $1/N < \mu$. Hence, for $n\geq N$, \quad $m\geq n$,
\begin{eqnarray*}
G(a_{n}, a_{m}, \lambda) &\geq & G(a_{n}, a_{m}, 1/n)\geq 1- 1/n > 1-\mu, \\
B(a_{n}, a_{m}, \lambda) &\leq & B(a_{n}, a_{m}, 1/n) \leq 1/n \leq \mu, \\
Y(a_{n}, a_{m}, \lambda) &\leq & Y(a_{n}, a_{m}, 1/n) \leq 1/n \leq \mu.
\end{eqnarray*}
Therefore, the sequence $(a_{n})$ is a Cauchy. We know that $F$ is complete. Then there exists $a\in F$ such that $a_{n} \rightarrow a$. Since $a_{k}\in O(a_{n}, \varepsilon_{n}^{*}, \lambda_{n}^{*})$ for
$k\geq n$, then we have $a\in \overline{O(a_{n}, \varepsilon_{n}^{*}, \lambda_{n}^{*})}$. Hence $a\in \overline{O(a_{n}, \varepsilon_{n}^{*}, \lambda_{n}^{*})}\subset O(a_{n-1}, \varepsilon_{n-1}^{*}, \lambda_{n-1}^{*})\cap \gamma_{n}$, \quad $\forall n$. Then, $\delta\cap (\cap_{n\in \mathbb{N}}\gamma_{n})\neq \emptyset$. Then $\cap_{n\in \mathbb{N}}\gamma_{n}$ is dense in $F$.
\end{proof}

\begin{Definition}
Let $(F, \mathcal{N}, \circ, \bullet)$ be a NMS. A collection $(\mathcal{D}_{n})$ $(n\in \mathbb{N})$ is said to have neutrosophic diameter zero (NDZ) if for each $0<\varepsilon<1$ and each $\lambda>0$, then there exists $N\in \mathbb{N}$ such that $G(a,b,\lambda)>1-\varepsilon$, $B(a,b,\lambda)<\varepsilon$ and $Y(a,b,\lambda)<\varepsilon$ for all $a,b\in \mathcal{D}_{N}$.
\end{Definition}

\begin{Theorem}
The NMS $(F, \mathcal{N}, \circ, \bullet)$ is complete if and only if every nested sequence $(\mathcal{D}_{n})_{n\in \mathbb{N}}$ of nonempty closed sets with NDZ have nonempty intersection.
\end{Theorem}

\begin{proof}
 Firstly consider the given condition is satisfied. We will show that $(F, \mathcal{N}, \circ, \bullet)$ is complete. Choose the Cauchy sequence $(a_{n})$ in $F$. If we define the $\mathcal{E}_{n}=\lbrace a_{k}: k \geq n \rbrace$ and $\mathcal{D}_{n}=\overline{\mathcal{E}_{n}}$, then we can say that $(\mathcal{D}_{n})$ has NDZ. For given $\zeta\in (0,1)$ and $\lambda >0$, we take $\varepsilon\in (0,1)$ such that $(1-\varepsilon)\bullet (1-\varepsilon)\bullet (1-\varepsilon)> 1-\zeta$ and $\varepsilon\bullet \varepsilon\bullet \varepsilon<\zeta$. Since the sequence $(a_{n})$ is Cauchy, then there exist $N \in \mathbb{N}$ such that $G(a_{n}, a_{m}, \frac{\lambda}{3})>1-\varepsilon$, $B(a_{n}, a_{m}, \frac{\lambda}{3})<\varepsilon$ and $Y(a_{n}, a_{m}, \frac{\lambda}{3})<\varepsilon$, \quad $(\forall m,n\geq N)$. Then, $G(a,b, \frac{\lambda}{3})>1-\varepsilon$, $B(a,b, \frac{\lambda}{3})<\varepsilon$ and $Y(a,b, \frac{\lambda}{3})<\varepsilon$, \quad $(\forall m,n\geq \mathcal{E_{N}})$.\\

  Choose $a,b\in \mathcal{D}_{N}$. There exist the sequences $(a_{n}^{*})$ and $(b_{n}^{*})$ such that $a_{n}^{*}\rightarrow a$ and $b_{n}^{*}\rightarrow b$. Thus, for sufficiently large $n$, $a_{n}^{*}\in O(a,\varepsilon, \frac{\lambda}{3})$ and $b_{n}^{*}\in O(b,\varepsilon, \frac{\lambda}{3})$. Now, we have
 \begin{eqnarray*}
 && G(a,b,\lambda)\geq G(a,a_{n}^{*},\frac{\lambda}{3})\circ G(a_{n}^{*}, b_{n}^{*}, \frac{\lambda}{3})\circ G(b_{n}^{*}, b, \frac{\lambda}{3})> (1-\varepsilon)\circ (1-\varepsilon)\circ (1-\varepsilon)>1-\zeta, \\
 &&  B(a,b,\lambda)\leq B(a,a_{n}^{*},\frac{\lambda}{3})\bullet B(a_{n}^{*}, b_{n}^{*}, \frac{\lambda}{3})\bullet B(b_{n}^{*}, b, \frac{\lambda}{3})< \varepsilon\bullet \varepsilon \bullet \varepsilon < \zeta,\\
 &&  Y(a,b,\lambda)\leq Y(a,a_{n}^{*},\frac{\lambda}{3})\bullet Y(a_{n}^{*}, b_{n}^{*}, \frac{\lambda}{3})\bullet Y(b_{n}^{*}, b, \frac{\lambda}{3})< \varepsilon\bullet \varepsilon \bullet \varepsilon < \zeta.
 \end{eqnarray*}
  From here, $G(a,b,\lambda)>1-\zeta$, $B(a,b,\lambda)<\zeta$ and $Y(a,b,\lambda)<\zeta$ \quad $(\forall a,b\in \mathcal{D}_{N})$. Therefore, $(\mathcal{D}_{N})$ has NDZ and so by the hypothesis $\cap_{n\in \mathbb{N}}\mathcal{D}_{n}$ is nonempty. Take $a\in \cap_{n\in \mathbb{N}}\mathcal{D}_{n}$. For $\varepsilon \in (0,1)$ and $\lambda >0$, then there exist $N_{1} \in \mathbb{N}$ such that $G(a_{n},a,\lambda)>1-\varepsilon$, $B(a_{n},a,\lambda)<\varepsilon$ and $Y(a_{n},a,\lambda)<\varepsilon$ \quad $(\forall n\geq N_{1})$. Therefore, for each $\lambda>0$, $G(a_{n}, a, \lambda)\rightarrow 1$, $B(a_{n}, a, \lambda)\rightarrow 0$ and $Y(a_{n}, a, \lambda)\rightarrow 0$ as $n\rightarrow \infty$. Hence, $a_{n}\rightarrow a$, that is
  $(F, \mathcal{N}, \circ, \bullet)$ is complete.\\

  Conversely, assume that $(F, \mathcal{N}, \circ, \bullet)$ is complete. Let's $(\mathcal{D}_{n})_{n\in \mathbb{N}}$  is nested sequence of nonempty closed sets with NDZ. For each $n\in \mathbb{N}$, take a point $a_{n}\in \mathcal{D}_{n}$. We will show that the sequence $(a_{n})$ is Cauchy. Since $(\mathcal{D}_{n})$ has NDZ, for $\lambda >0$ and $0<\varepsilon<1$, then there exist $N \in \mathbb{N}$ such that $G(a,b,\lambda)>1-\varepsilon$, $B(a,b,\lambda)<\varepsilon$ and $Y(a,b,\lambda)<\varepsilon$ \quad $(\forall a,b\in \mathcal{D}_{N})$. Since the sequence $(D_{n})$ is nested, then $G(a_{n}, a_{m}, \lambda)>1-\varepsilon$, $B(a_{n}, a_{m}, \lambda)<\varepsilon$ and $Y(a_{n}, a_{m}, \lambda)<\varepsilon$ \quad $(\forall m,n \geq N)$. Hence the sequence $(a_{n})$ is Cauchy. Since $(F, \mathcal{N}, \circ, \bullet)$ is complete, then $a_{n}\rightarrow a$ for some $a\in F$. Therefore, $a\in \overline{\mathcal{D}_{n}}=\mathcal{D}_{n}$ for every $n$, and so $a\in \cap_{n\in \mathbb{N}}\mathcal{D}_{n}$.
\end{proof}

\begin{Theorem}
Every separable NMS is second countable.
\end{Theorem}

\begin{proof}
  Give the separable NMS $(F, \mathcal{N}, \circ, \bullet)$. Let $A=\lbrace a_{n}: n\in \mathbb{N} \rbrace$ be a countable dense subset of $F$. Establish the family $\textbf{O}=\lbrace O(a_{k}, 1/m, 1/m):  k,m\in \mathbb{N} \rbrace$. 
It can be easily seen, $\textbf{O}$ is countable. We will show that $\textbf{O}$ is base for the family of all OSs in $F$. Let $\gamma$ be any OS in $F$, \quad $a\in \gamma$. Then there exist $\lambda >0$, \quad $0<\varepsilon <1$ such that $O(a,\varepsilon, \lambda)\subset \gamma$. Since $0<\varepsilon <1$, we can choose a $0<\zeta<1$ such that $(1-\zeta)\circ(1-\zeta)>1-\varepsilon$ and $\zeta\bullet \zeta<\varepsilon$. Take $t\in \mathbb{N}$ such that $1/t<\min\lbrace \zeta, \lambda /2 \rbrace$. 
Since it is known that $A$ is dense in $F$, there exist $a_{k}\in A$ such that $a_{k}\in O(a, 1/t, 1/t)$. If $b\in O(a_{k}, 1/t, 1/t)$, we have
  \begin{eqnarray*}
  G(a,b,\lambda)\geq G(a,a_{k}, \frac{\lambda}{2})\circ G(b,a_{k}, \frac{\lambda}{2}) \geq G(a,a_{k}, \frac{1}{t}) \circ G(b,a_{k}, \frac{1}{t}) \\
  \geq (1- \frac{1}{t})\circ (1-\frac{1}{t})\geq (1-\zeta) \circ (1-\zeta) > 1-\varepsilon, \\
   B(a,b,\lambda)\leq B(a,a_{k}, \frac{\lambda}{2})\bullet B(b,a_{k}, \frac{\lambda}{2}) \leq B(a,a_{k}, \frac{1}{t}) \bullet B(b,a_{k}, \frac{1}{t}) \leq \frac{1}{t} \bullet \frac{1}{t} \leq \zeta \bullet \zeta < \varepsilon,\\
   Y(a,b,\lambda)\leq Y(a,a_{k}, \frac{\lambda}{2})\bullet Y(b,a_{k}, \frac{\lambda}{2}) \leq Y(a,a_{k}, \frac{1}{t}) \bullet Y(b,a_{k}, \frac{1}{t}) \leq \frac{1}{t} \bullet \frac{1}{t} \leq \zeta \bullet \zeta < \varepsilon,
  \end{eqnarray*}
  Then, $b\in O(a,\varepsilon, \lambda)\subset \gamma$ and so $\textbf{O}$ is a base.
\end{proof}

Note that the second countability implies separability and the second countability is inheritable property. Then, we can say that every subspace of a separable NMS is separable.

\begin{Definition}
Let $F$ be any nonempty set and $(H, \mathcal{N}, \circ, \bullet)$ be a NMS. The sequence of functions $(f_{n}):F \rightarrow G$ is called converge uniformly to a function $f:F\rightarrow G$, if given $\lambda >0$, \quad $\varepsilon \in (0,1)$, then there exists $N\in \mathbb{N}$ such that $G(f_{n}(a), f(a), \lambda)>1-\varepsilon$, $B(f_{n}(a), f(a), \lambda)<\varepsilon$,
$Y(f_{n}(a), f(a), \lambda)<\varepsilon$ $\forall n\geq N$ and $\forall a\in F$.
\end{Definition}

Now, we will give Uniform Convergence Theorem for NMS:

\begin{Theorem}
Let $f_{n}:F\rightarrow H$ be a sequence of continuous functions from a topological space $F$ to a NMS $(H, \mathcal{N}, \circ, \bullet)$. If $(f_{n})$ converges uniformly to $f:F\rightarrow H$, then $f$ is continuous.
\end{Theorem}

\begin{proof}
  Take $\delta$ be OS of $H$ and let $a_{0}\in f^{-1}(\delta)$. Since $\delta$ is open, then there exist $\lambda >0$, \quad $\varepsilon\in (0,1)$ such that $O(f(a_{0}), \varepsilon, \lambda)\subset \delta$. Since $\varepsilon\in (0,1)$, we take a $\zeta\in (0,1)$ such that $(1-\zeta)\circ (1-\zeta)\circ (1-\zeta)>1-\varepsilon$ and $\zeta\bullet \zeta\bullet \zeta< \varepsilon$. Since $(f_{n})$ converges uniformly to $f$, then, for $\lambda >0$, \quad $\zeta\in (0,1)$, there exists $N \in \mathbb{N}$ such that $G(f_{n}(a), f(a), \frac{\lambda}{3})>1-\zeta$, $B(f_{n}(a), f(a), \frac{\lambda}{3})<\zeta$ and $Y(f_{n}(a), f(a), \frac{\lambda}{3})<\zeta$ $\forall n\geq N$ and $\forall a\in F$. Since $f_{n}$ continuous $\forall n\in \mathbb{N}$, then there exist a neighborhood $\gamma$ of $a_{0}$ such that $f_{n}(\gamma)\subset O(f_{n}(a_{0}), \zeta, \frac{\lambda}{3})$. Hence $G(f_{n}(a), f_{n}(a_{0}), \frac{\lambda}{3})>1-\zeta$, $B(f_{n}(a), f_{n}(a_{0}), \frac{\lambda}{3})<\zeta$ and $Y(f_{n}(a), f_{n}(a_{0}), \frac{\lambda}{3})<\zeta$ for all $a\in \gamma$. Now
  \begin{eqnarray*}
  G(f(a), f(a_{0}), \lambda) \geq G(f(a), f_{n}(a), \frac{\lambda}{3})\circ G(f_{n}(a), f_{n}(a_{0}), \frac{\lambda}{3})\circ G(f_{n}(a_{0}), f(a_{0}), \frac{\lambda}{3})\\
  \geq (1-\zeta)\circ (1-\zeta)\circ (1-\zeta)>1-\varepsilon,\\
  B(f(a), f(a_{0}), \lambda) \leq B(f(a), f_{n}(a), \frac{\lambda}{3}) \bullet B(f_{n}(a), f_{n}(a_{0}), \frac{\lambda}{3})\bullet B(f_{n}(a_{0}), f(a_{0}), \frac{\lambda}{3})\leq \zeta \bullet \zeta \bullet \zeta < \varepsilon, \\
  Y(f(a), f(a_{0}), \lambda) \leq Y(f(a), f_{n}(a), \frac{\lambda}{3}) \bullet Y(f_{n}(a), f_{n}(a_{0}), \frac{\lambda}{3})\bullet Y(f_{n}(a_{0}), f(a_{0}), \frac{\lambda}{3})\leq \zeta \bullet \zeta \bullet \zeta < \varepsilon.
  \end{eqnarray*}
  Therefore, $f(a)\in O(f(a_{0}), \varepsilon, \lambda)\subset \delta$ for all $a\in \gamma$. Hence $f(\gamma)\subset \delta$ and so $f$ is continuous.
\end{proof}

\section{Conclusion}
The aim of this study is to define a neutrosophic metric spaces and examine some properties. The structural characteristic properties of NMSs such as open ball, open set, Hausdorffness, compactness, completeness, nowhere dense in NMS have been established. Analogues of Baire Category Theorem and Uniform Convergence Theorem  are given for NMS.\\

\end{document}